\documentclass[11pt]{article}
\usepackage{mathrsfs}

\input{amssymb.sty}
\usepackage{pb-diagram}
\begin{document}

\title{\bf Perturbation of  Browder Spectrum of Upper-triangular Operator Matrices\footnote{This work is
supported by  the NSF of China (Grant Nos. 11301077,11301124,
11226113, 11171066 and 10771191 ) , the FED of Fujian
Province(Grant No. JA12074), the NSF of Fujian Province (Grant
No.2012J05003).}}

\author{{Shifang Zhang$^{1}$\footnote{Corresponding author:
E-mail: shifangzhangfj@163.com}, \, Huaijie Zhong$^{1}$, \, Lin
 Zhang $^2$}
\\\\$^1$\small\it School of Mathematics and Computer Science, Fujian Normal University, Fuzhou 350007, P. R. China
\\\\$^2$\small\it  Institute of Mathematics, Hangzhou Dianzi University, Hangzhou 310018, P. R. China}

\date{} \maketitle
\begin{center}
\begin{minipage}{140mm}
{\bf Abstract} { Let
$M_{C}=\left(\begin{array}{cc}A&C\\0&B\\\end{array} \right)$ is a
2-by-2 upper triangular operator matrix acting on the Banach space
$X\oplus Y$ or Hilbert space $H\oplus K$. For the most import
spectra such as spectrum, essential spectrum and Weyl spectrum,
the characterizations of the perturbations of $M_{C}$ on Banach
space have been presented. However, the characterization for the
Browder spectrum on the Banach space is still unknown. The goal of
this paper is to present some necessary and sufficient conditions
for $M_C$ to be Browder for some $C\in B(Y,X)$ by using an
alternative approach based on matrix representation of operators
and the ghost index theorem. Moreover, in the Hilbert space case
the characterizations of the Fredholm and invertible perturbations
of Browder spectrum are also given. }

\vspace{4mm} {\bf Keywords}: \,{operator matrices, Browder spectrum,
left semi-Browder spectrum.}
\end{minipage}
\end{center}

\vskip 0.2in
%
%
 \large\section{Introduction}

It is well-known that if $H$ is a Hilbert space, $T$ is a bounded
linear operator on $H,$ and $H_1$ is an invariant closed subspace of
$T$, then $T$ is of the form
$$T=\left(
       \begin{array}{cc}
        *&*\\
         0&*\\
       \end{array}
     \right):H_1\oplus H_1^{\perp}\rightarrow H_1\oplus H_1^{\perp},$$
which motivated the interest in $2 \times 2$ upper-triangular
operator matrices (see [1-9, 11-20]). The problems related to the
perturbation of spectra of  2-by-2 upper triangular operator
matrices were first studied  by H. K. Du and J. Pan in [5], and
they  given the characterization of $M_C$ to be invertible for
some $C\in B(K,H)$ in the Hilbert spece case. Some years later, in
[13] J. K. Han et al. generalized  the above result to the  Banach
space case.  In 2002, D. S. Djordjevi\'{c} in [5] future studied
this field for some other spectra and  presented the necessary and
sufficient conditions for $M_C$ to be a Fredholm and Weyl operator
for some $C\in B(Y,X)$, respectively. Also in [5], some sufficient
conditions for $M_C$ to be a Browder operator were given.
Recently, Cao [1] gave the necessary and sufficient condition for
$M_C$ to be a Browder operator in the Hilbert space case. Till
now, although for the most import and widely used spectra such as
spectrum, essential spectrum and Weyl spectrum, the
characterizations of the perturbations of upper triangular
operator matrices on Banach space have been given, the
characterization for the Browder spectrum is still unknown. The
goal of  this paper is to present a necessary and sufficient
condition  for $M_C$ to be Browder for some $C\in B(Y,X)$ by using
an alternative approach based on matrix representation of
operators and the ghost index theorem,  and our result is an
extension of the main result in [1] and the corresponding results
from [5, 22]. Moreover, the characterizations of the Fredholm and
invertible perturbations of Browder spectrum are also given.

Throughout this paper, let $H$ and $K$ be the complex infinite
dimensional separable Hilbert spaces, $X$ and $Y$ be the complex
infinite dimensional Banach spaces, and $B(X, Y)$ be the set of
all bounded linear operators from $X$ into $Y.$ For simplicity, we
write $B(X, X)$ as $B(X)$. For $T \in B(X)$, $T$ is called Drazin
invertible if there exists an operator $T^D \in B(X)$ such that
$$TT^D=T^D T ,\quad \quad T^DTT^D=T^D ,\quad \quad T^{k+1}T^D=T^{k},$$
\noindent for some nonnegative integer $k$([5]). $T\in B(X, Y)$ is
called regular( also called relatively regular in [5]) if there is
an operator $B\in B(Y, X)$ such that  $T=TBT$([4]). It is
well-known that $T$ is regular if and only if $R(T)$ and $N(T)$,
respectively, are closed and complemented subspaces of $Y$ and
$X$.

For $T\in B(X, Y)$, let $R(T)$ and $N(T)$ be the range and kernel
of $T$, respectively, and denote $\alpha (T)=\dim N(T)$ and $\beta
(T)=\dim Y / R(T)$. If $T\in B(X)$, the ascent, denoted by
$asc(T)$, and the descent, denoted by $des(T)$, of $T$ are the
non-negative integers $k$ defined  respectively by
$$asc(T)=inf\{k:N(T^{k})=N(T^{k+1})\},\,\,\,\,\,\,des(T)=inf\{k:R(T^{k})=R(T^{k+1})\}.$$
If such $k$ does not exist, then $asc(T)=\infty$, respectively
$des(T)=\infty$. If the ascent and the descent of $T$ are finite,
then they are equal (see [5]). For $T\in B(X)$, if $R(T)$ is
closed and $\alpha (T)<\infty$, then $T$ is said to be  upper
semi-Fredholm, if $\beta (T)<\infty$, then $T$ is said to be
semi-Fredholm. If $T\in B(X)$ is either upper or lower
semi-Fredholm , then $T$  is said to be semi-Fredholm. For a
semi-Fredholm operator $T$, its index ind $(T)$ is defined as ind
$(T)=\alpha(T )-\beta(T).$ The sets of invertible operators, left
invertible operators and right invertible operators on $X$,
respectively, are denoted by
 $G(X), G_l(X)$ and $G_r(X).$

The sets of all Fredholm operators, Weyl operators, left
semi-Fredholm operators and  right semi-Fredholm operators on $X$
are defined, respectively, by
\begin{eqnarray*}
& & \Phi(X):=\{T \in  B(X):\alpha (T)<\infty \makebox{ and }\beta (T)< \infty \};\\
&& \Phi_{0}(X):=\{T \in \Phi(X): \makebox {ind} (T)= 0 \};\\
& & \Phi_{l}(X):=\{T \in  B(X):
 R(T)\makebox{ is a closed and complemented subspace of }X  \makebox{ and }\,\,\alpha (T)<\infty\};\\
& &\Phi_r(X):=\{T \in  B(X):
 N(T)\makebox{ is a closed and complemented subspace of }X  \makebox{ and }\,\,\beta (T)<\infty \}.
\end{eqnarray*}

The sets of all Browder operators,  left semi-Browder operators and
right semi-Browder operators on $X$ are defined, respectively, by
\begin{eqnarray*}
& & \Phi_{b}(X):=\{T \in  \Phi(X): asc(T)= des(T)<\infty \};\,\,\,\,\,\,\,\,\,\,\,\,\,\,\,\,\,\,\,\,\,\,\,\,\,\,\,\,\,\,\\
& &  \Phi_{lb}(X):=\{T \in  \Phi_l(X): asc(T)<\infty \};\\
& &\Phi_{rb}(X):=\{T \in  \Phi_r(X): des(T)<\infty \}.
\end{eqnarray*}
The spectrum, left spectrum, right spectrum, essential spectrum,
Weyl spectrum, left semi-Fredholm spectrum, right semi-Fredholm
spectrum, Browder spectrum, left semi-Browder spectrum, right
semi-Browder spectrum of an operator $T\in B(H)$ are respectively
denoted by $\sigma(T)$, $\sigma_{l}(T)$, $\sigma_{r}(T)$,
$\sigma_{e}(T)$, $\sigma_{w}(T)$,$\sigma_{le}(T)$,
$\sigma_{re}(T)$, $\sigma_{b}(T)$, $\sigma_{lb}(T)$ and
$\sigma_{rb}(T)$.

\vspace{4mm}

One of the main results is : $M_{C}$ is a Browder operator for
some $C\in B(Y,X)$ if and only if the following statements are
satisfied
\par (a) $A\in B(X)$ is a left semi-Browder operator;
\par (b) $B\in B(Y)$ is a right semi-Browder operator;
\par (c) $ N(A)\times N(B)\cong X / R(A)\times Y / R(B).$

\noindent Using the above result, we characterize completely the set
$\bigcap_{C\in B(Y,\,X)}\sigma_{b}(M_{C})$, which is the extension
 of the corresponding results from [1, 5, 22], that is,
$$\bigcap_{C\in B(Y,\,X)}\sigma_{b}(M_{C})=\sigma_{lb} (A)\cup \sigma_{rb}(B)\cup\{\lambda\in {\mathbb{C}}:
N(A-\lambda)\times N(B-\lambda)\not\cong X / R(A-\lambda)\times Y /
R(B-\lambda)\}.$$
Moreover, for any given $A\in B(H)$ and $ B\in
B(K)$, it is proved that

\begin{eqnarray*} \bigcap_{C\in G(K,\,H)}\sigma_{b}(M_{C})&=& \bigcap_{C\in \Phi(K,\,H)} \sigma_{b}(M_{C})\\
&=& \sigma_{lb} (A)\cup \sigma_{rb}(B)\cup\{\lambda\in {\mathbb{C}}:
\alpha(A-\lambda)+ \alpha(B-\lambda)\not= \beta(A-\lambda)+
\beta(B-\lambda)\}.\end{eqnarray*}

%
%
%
\vskip 0.2in
 \large\section{Main Results and Proofs}

\par We begin with some lemmas, which are useful for the proofs of  the main results.

\vskip 0.2in \noindent {\bf Lemma 2.1.} ([23] or [2]) For any
given $A\in B(X)$ and $ B\in B(Y)$, if $M_{C}$ is Drazin
invertible for some $C\in B(Y,\,X)$, then
 \par (1)~$des(B)<\infty$ and $asc(A)<\infty;$
 \par (2)~$des(A^*)<\infty$ and $asc(B^*)<\infty$.

\vskip 0.2in \noindent {\bf Lemma 2.2.} ([5] Lemma 2.3) Let $M$
and $N$ be finite dimensional spaces. If $\dim M=\dim N$ and   $
X\times M\cong Y\times N$, then  $ X\cong Y.$

\vskip 0.2in \noindent {\bf Lemma 2.3.} ([13]) If $T\in B(X, Y)$,
$S\in B(Y, Z)$ and  $ST\in B(X, Z)$ are  regular, then

$$ N(T)\times  N(S)\times Z/ R(ST)\cong  N(ST)\times Y/ R(T)\times Z / R(S).$$

Some well-known facts are collected in the following lemma.

\medskip\noindent {\bf Lemma 2.4.} ([24])
 For any given $A\in B(X)$, $B\in B(Y)$ and $C\in B(Y,X)$, we have:

\par (1)~ if any two of operators $A, B$ and $M_C$ are invertible (resp., Fredholm, Weyl, Browder,Drazin invertible ), then so is the third;
\par (2) if  $A$ is Browder, then $B$ is left semi-Browder if and only if so is $M_C$ ;
\par (3)~ if  $B$ is Browder, then $A$ is right semi-Browder if and only if  so is $M_C$.

\vskip 0.2in \noindent {\bf Lemma 2.5}([5]). Let $A\in B(X)$ and $
B\in B(Y)$ be given and consider the statements:
\par (i)~ $M_{C}$ is a Weyl operator for some $C\in B(Y,X)$
\par (ii)~ $A\in B(X)$ and $B\in B(Y)$ satisfy the following conditions:
\par (a)~ $A$ is a left semi-Fredholm operator;
\par (b)~ $B$ is a right semi-Fredholm operator;
\par (c)~ $ N(A)\times N(B)\cong X / R(A)\times Y / R(B).$
\par Then (i)~ $\Longleftrightarrow$ (ii)

\vskip 0.2in \noindent {\bf Lemma 2.6} ([13]).  Let $A\in B(X)$
and $ B\in B(Y)$ be given and consider the statements:
\par (i) $M_{C}$ is invertible for some $C\in B(Y,X)$
\par (ii) $A\in B(X)$ and $B\in B(Y)$ satisfy the following conditions:
\par (a) $A$ is left invertible;
\par (b) $B$ is right invertible;
\par (c) $ N(B)\cong X / R(A).$
\par Then (i) $\Leftrightarrow$ (ii)

\vskip 0.2in \noindent {\bf Lemma 2.7}${\label{2}}$. For $T\in
B(X)$, $T$ is left semi-Browder if and only if $T$ can be decomposed
into the following form with respect to space decomposition
$X=X_1\oplus X_2$
$$T=\left(\begin{array}{cc}T_1&T_{12}\\0&T_2\end{array}\right),$$
\noindent where $\dim(X_1)<\infty$, $T_1$ is nilpotent, and $T_2$ is
 left invertible.

 \vskip 0.2in \noindent {\bf Proof.}
Necessity. Suppose that  $T$ is  left semi-Browder. Then we can
assume $p=asc(T)<\infty$. Let $X_1= N(T^p)$. Then $\dim X_1<\infty$
since $T$ is left semi-Fredholm. Therefore, there exists a closed
subspace $X_2$ of $X$ such that $X_1\oplus X_2=X$. With respect to
the space decomposition $X=X_1\oplus X_2$, we have
 $$T=\left( \begin{array}{cc}  T_1&T_{12}\\ 0&T_2\\ \end{array}
 \right):X_1\oplus X_2 \rightarrow  X_1\oplus X_2.$$
Obviously, $T_1$ is nilpotent and hence $T_1\in\Phi_{b}(X_1)$.
Moreover, it follows from Lemma 2.4 that $T_2\in\Phi_{lb}(X_2).$
Next, we shall  show $T_2$ is injective. If not, there exists some
$0\not=x\in X_2$ such that $T_2x=0.$ Put $z=\left(\begin{array}{c}
   0\\x \\ \end{array}\right)$. Then $T^{p+1}\left(\begin{array}{c}
   0\\x \\ \end{array}\right)=T^{p}\left(\begin{array}{c}
   T_{12}x\\0 \\ \end{array}\right)=0$, and so $0\not=\left(\begin{array}{c}
   0\\x \\ \end{array}\right)\in N(T^{p+1})$. This is in
contradiction with the assumption that $asc(T)=p<\infty$, thus,
$T_2$ is left invertible.

 Sufficiency. It is evident.

\vskip 0.2in \noindent {\bf Lemma 2.8}${\label{2}}$. For $T\in
B(X)$, $T$ is right semi-Browder if and only if $T$ can  be
decomposed into the following  form with respect to some space
decomposition $X=X_1\oplus X_2$

$$T=\left(\begin{array}{cc}T_1&T_{12}\\0&T_2\end{array}\right),$$

\noindent where $\dim(X_2)<\infty$, $T_1$ is  right invertible, and
$T_2$ is nilpotent.

 \vskip 0.2in \noindent {\bf Proof.}
 Necessity. Suppose  $T$ is  right semi-Browder.
 Then we can assume $p=des(T)<\infty$. Now  put $X_1= R(T^p)$. Then $\dim
 X /X_1<\infty$ since $T^p$ is a right semi-Browder operator. Therefore, there
 exists  a closed finite dimensional subspace $X_2$ of $X$ such that
 $X_1\oplus X_2=X$. Corresponding the space decomposition $X=X_1\oplus X_2$,
 we have
   $$T=\left( \begin{array}{cc}   T_1&T_{12}\\ 0&T_2\\ \end{array}
   \right):X_1\oplus X_2 \rightarrow  X_1\oplus X_2.$$
Obviously, $T_1$ is surjective and $T_2^P=0$. Moreover, noting that
$\dim X_2<\infty$ implies $T_2\in\Phi_{b}(X_2)$, it follows from
Lemma ~2.4 that $T_1\in\Phi_{rb}(X_2)$, and so $T_1$ is right
invertible.

 Sufficiency. It is evident.

\vskip 0.2in

One of our main results is as follows

\vskip 0.2in \noindent {\bf Theorem 2.9}. Let $A\in B(X)$ and $ B\in
B(Y)$ be given and consider the statements:
\par (i) $M_{C}$ is a Browder operator for some $C\in B(Y,X)$
\par (ii) $A\in B(X)$ and $B\in B(Y)$ satisfy the following conditions:
\par (a) $A$ is a left semi-Browder  operator;
\par (b) $B$ is a right semi-Browder  operator;
\par (c) $ N(A)\times N(B)\cong X / R(A)\times Y / R(B).$
\par Then (i) $\Leftrightarrow$ (ii).

\vskip 0.2in \noindent {\bf Proof.} Necessity. Noting the fact that
an operator $T$ is a Browder operator if and only if $T$ is a Drazin
invertible Weyl operator, the implication  follows from Lemmas 2.1
and 2.5 immediately.

Sufficiency. Since $A$ is a left semi-Browder operator, from Lemma
2.7, we have

$$A=\left( \begin{array}{cc} A_1&A_{12}\\0&A_2\\\end{array}
 \right):X_1\oplus X_2 \rightarrow  X_1\oplus X_2,$$
where $A_1$ is nilpotent, $A_2$ is left invertible and $\dim
(X_1)<\infty.$ Note that $$A=\left(\begin{array}{cc} I&0\\ 0&A_2\\
\end{array}\right)\left(\begin{array}{cc}   A_1&A_{12}\\ 0&I\\
\end{array}\right).$$ Applying Lemma 2.3, we have

$$N(\left( \begin{array}{cc}   I&0\\ 0&A_2\\ \end{array}\right))\times
N(\left(\begin{array}{cc}   A_1&A_{12}\\ 0&I\\ \end{array}\right))
\times X / R(A)\cong X /R\left( \begin{array}{cc} I&0\\ 0&A_2\\
\end{array}\right)\times X /R\left( \begin{array}{cc} A_1&A_{12}\\
0&I\\\end{array}\right)  \times N(A).$$
Since $\left(
\begin{array}{cc}   I&0\\ 0&A_2\\ \end{array}\right)$ is
injective,

$$ N(\left(\begin{array}{cc}   A_1&A_{12}\\ 0&I\\\end{array}\right))
\times X / R(A)\cong X /R\left( \begin{array}{cc} I&0\\ 0&A_2\\
\end{array}\right)\times X /R\left(\begin{array}{cc} A_1&A_{12}\\
 0&I\\ \end{array}\right) \times N(A).$$
It follows from $\dim (X_1)<\infty$ that $A_1\in\Phi_{b}(X_1).$
By Lemma 2.4, we obtain $\left(\begin{array}{cc}   A_1&A_{12}\\
0&I\\ \end{array}\right)\in \Phi_{b}(X_1\oplus X_2)$, which implies

$$\alpha(\left(\begin{array}{cc}   A_1&A_{12}\\ 0&I\\
\end{array}\right))= \beta(\left(\begin{array}{cc}
 A_1&A_{12}\\ 0&I\\ \end{array}\right))<\infty.$$
Using Lemma 2.2, we have

$$ X / R(A)\cong X /R\left( \begin{array}{cc} I&0\\ 0&A_2\\
\end{array}\right)\times N(A),$$
and hence \begin{equation} X / R(A) \cong X_2 /R(A_2) \times
N(A).\end{equation}

\noindent Meanwhile, since $B\in B(Y)$ is a right semi-Browder
operator, from Lemma 2.8 we have

$$B=\left( \begin{array}{cc}   B_1&B_{12}\\ 0&B_2\\ \end{array}
\right):Y_1\oplus Y_2 \rightarrow  Y_1\oplus Y_2,$$

\noindent where $B_2$ is nilpotent, $B_1$ is right invertible and
$\dim (Y_2)<\infty.$ Similar to the preceding arguments, we have

$$B=\left( \begin{array}{cc} I&B_{12}\\ 0&B_2\\ \end{array}\right)
\left(\begin{array}{cc}   B_1&0\\ 0&I\\ \end{array}\right),$$

\noindent and

$$N(\left( \begin{array}{cc} I&B_{12}\\ 0&B_2\\\end{array}\right))
\times N(\left(\begin{array}{cc}B_1&0\\ 0&I\\
\end{array}\right))\times Y / R(B)\cong Y /R\left( \begin{array}{cc}
I&B_{12}\\ 0&B_2\\ \end{array}\right)
 \times Y/R\left( \begin{array}{cc} B_1&0\\ 0&I\\ \end{array}\right)
  \times N(B).$$
Since $\left( \begin{array}{cc} B_{1}&0\\ 0&I\\
\end{array}\right)$ is surjective,
 $$N(\left( \begin{array}{cc} I&B_{12}\\ 0&B_2\\ \end{array}\right))\times
 N(\left(\begin{array}{cc} B_1&0\\ 0&I\\ \end{array}\right))\times Y / R(B)
 \cong Y /R\left( \begin{array}{cc} I&B_{12}\\ 0&B_2\\ \end{array}\right)
  \times N(B).$$

\noindent  It follows from $\dim (Y_2)<\infty$ that
$B_2\in\Phi_{b}(Y_2).$ By Lemma 2.4, we have

 $$\left( \begin{array}{cc}I&B_{12}\\ 0&B_2\\ \end{array}\right)\in\Phi_{b}(Y),$$

\noindent which implies $$\alpha(\left(\begin{array}{cc}I&B_{12}\\
0&B_2\\ \end{array}\right))=\beta(\left(\begin{array}{cc} I&B_{12}\\
0&B_2\\ \end{array}\right))<\infty.$$
 Then it follows from Lemma 2.2 that
\begin{equation} N(B_1)\times Y / R(B) \cong  N(B).\end{equation}

\noindent Combining (1) and (2) with the assumption $ N(A)\times
N(B)\cong X / R(A) \times Y / R(B)$, it is easy to show that
$$ N(A)\times N(B_1)\times Y / R(B)\cong X_2 /R(A_2) \times N(A) \times
Y / R(B).$$

\noindent Moreover, since  $\alpha(A)$ and $\beta (B)$ are finite,
it results from Lemma 2.2 that \begin{equation}   N(B_1)\cong X_2
/R(A_2).\end{equation}   Observe that $B_1$ is right invertible and
$A_2$ is left invertible. It follows from Lemma 2.6 that there
exists some $C_0\in B(Y_1, X_2)$ such that
$$\left(\begin{array}{cc}   A_2&C_0\\ 0&B_1\\ \end{array}
  \right) \in G(X_2\oplus Y_1).$$
Define $C\in B( Y, X)$ as $$C=\left( \begin{array}{cc}   0&0\\ C_0&0\\
\end{array} \right):Y_1\oplus Y_2 \rightarrow X_1\oplus X_2.$$
Then $M_C$ can be written as the following form:

$$M_C=\left( \begin{array}{cccc} A_1&B_{12}&0&0\\0&A_2&C_0&0\\0&0&B_1&
B_{12}\\0&0&0&B_2\\ \end{array}\right):X_1\oplus X_2\oplus Y_1\oplus
Y_2 \rightarrow  X_1\oplus X_2\oplus Y_1\oplus Y_2.$$

\noindent We conclude from Lemma 2.4 that $M_C=\left(
\begin{array}{cc} A&C\\0&B\\\end{array}  \right)\in\Phi_{b}(X\oplus
Y).$ This completes the proof.

\vskip 0.2in The next corollary is immediate from Theorem 2.9, which
characterizes the set $\bigcap_{C\in B(Y,\,X)}\sigma_{b}(M_{C}).$

\vskip 0.2in \noindent {\bf Corollary 2.10.} Let $A\in B(X)$ and $
B\in B(Y)$. Then we have

$$\bigcap_{C\in B(Y,\,X)}\sigma_{b}(M_{C})=\sigma_{lb}(A)\cup
\sigma_{rb}(B)\cup\{\lambda\in {\mathbb{C}}:N(A-\lambda)\times
N(B-\lambda)\not\cong X /R(A-\lambda)\times Y / R(B-\lambda)\}.$$

As a particular case of Corollary 2.10, we can rewrite the main
result in [1]

\vskip 0.2in \noindent {\bf Corollary 2.11.}  Let $H$ and $K$ be the
complex infinite dimensional separable Hilbert space, $A\in B(H)$
and $ B\in B(K)$. Then
$$\bigcap_{C\in B(K,\,H)}\sigma_{b}(M_{C})=\sigma_{lb}(A)\cup
\sigma_{rb}(B)\cup\{\lambda\in {\mathbb{C}}:\alpha(A-\lambda)+
\alpha(B-\lambda)\not= \beta(A-\lambda)+ \beta(B-\lambda)\}.$$

 \vskip 0.2in In [25], we have studied the sets $\bigcap_{C\in
G(K,\,H)}\sigma_{*}(M_{C})$ and
$\bigcap_{C\in\Phi(K,\,H)}\sigma_{*}(M_{C}),$ where
$\sigma_{*}(\cdot)$ can be equal to spectrum, essential spectrum,
Weyl spectrum and so on. Similarly,  the sets  $\bigcap_{C\in
G(K,\,H)}\sigma_{b}(M_{C})$ and $\bigcap_{C\in
\Phi(K,\,H)}\sigma_{b}(M_{C})$ are studied in the next theorem.
Firstly, we give a well-known result which is needed in the proof
of the next theorem.

\vskip 0.2in \noindent {\bf Lemma 2.12.} [18 Corollary 2.2]  Let
$A\in B(H)$ and $ B\in B(K)$. If $R(B)$ is closed, then $M_C$ is
left invertible if and only if both $A$ and $C_1$ are left
invertible, where $C = \left(
       \begin{array}{cc}
        C_1&C_2\\
         C_3&C_4\\
       \end{array}
     \right):
       N(B)\oplus
         N(B)^\perp
      \rightarrow
     R(A)^\perp \oplus R(A)$.

\vskip 0.2in
Duality, we have
\vskip 0.2in

 \noindent {\bf Corollary 2.13.}  Let $A\in B(H)$ and $ B\in
B(K)$. If $R(A)$ is closed, then $M_C$ is right invertible if and
only if both $B$ and $C_1$ are right invertible, where $C = \left(
       \begin{array}{cc}
        C_1&C_2\\
         C_3&C_4\\
       \end{array}
     \right):
       N(B)\oplus
         N(B)^\perp
      \rightarrow
     R(A)^\perp \oplus R(A)$.

\vskip 0.2in Combing Lemma 2.12 and Corollary 2.13, we have
 \vskip 0.2in

\noindent {\bf Corollary 2.14.}  Let $A\in B(H)$ and $ B\in B(K)$.
Then $M_C$ is invertible if and only if  $A$ is left invertible, $B$
is right invertible  and $C_1$ is  invertible, where $C = \left(
       \begin{array}{cc}
        C_1&C_2\\
         C_3&C_4\\
       \end{array}
     \right):
       N(B)\oplus
         N(B)^\perp
      \rightarrow
     R(A)^\perp \oplus R(A)$.

\vskip 0.1in We are now ready to give  the last  main result.

 \vskip
0.2in \noindent {\bf Theorem 2.15.} For any given $A\in B(H)$ and $
B\in B(K)$ we have

\begin{equation}\bigcap_{C\in G(K,\,H)}\sigma_{b}(M_{C})=\bigcap_{C\in
\Phi(K,\,H)}\sigma_{b}(M_{C})=SPR,
\end{equation}
where $SPR = \sigma_{lb} (A)\cup \sigma_{rb}(B)\cup\{\lambda\in
{\mathbb{C}}: \alpha(A-\lambda)+ \alpha(B-\lambda)\not=
\beta(A-\lambda)+ \beta(B-\lambda)\}.$ \vskip 0.2in \noindent {\bf
Proof.} Noting the fact
$$ \bigcap_{C\in G(K,\,H)}\sigma_{b}(M_{C})\supseteq\bigcap_{C\in
\Phi(K,\,H)}\sigma_{b}(M_{C})\supseteq\bigcap_{C\in
B(K,\,H)}\sigma_{b}(M_{C}),$$ it follows from Corollary 2.11 that in
order to prove (4), we only need to prove
$$\bigcap_{C\in G(K,\,H)}\sigma_{b}(M_{C})\subseteq\sigma_{lb} (A)\cup
\sigma_{rb}(B)\cup\{\lambda\in {\mathbb{C}}: \alpha(A-\lambda)+
\alpha(B-\lambda)\not= \beta(A-\lambda)+
\beta(B-\lambda)\}.$$
 For this, it suffices to show that, if
 $$0\not \in \sigma_{lb} (A)\cup
\sigma_{rb}(B)\cup\{\lambda\in {\mathbb{C}}: \alpha(A-\lambda)+
\alpha(B-\lambda)\not= \beta(A-\lambda)+ \beta(B-\lambda)\},$$ then
$$0\not\in \bigcap_{C\in G(K,\,H)}\sigma_{b}(M_{C}).$$ Now suppose
that $A\in \Phi_{lb}(H)$ and $B\in \Phi_{rb}(K)$ with $\alpha(A)+
\alpha(B)= \beta(A)+ \beta(B).$ We shall prove that there exists
some $C\in G(K,H)$ such that $M_C \in \Phi_{b}(H\oplus K).$

From Lemmas 2.7 and 2.8, we have
$$A=\left( \begin{array}{cc} A_1&A_{12}\\ 0&A_2\\ \end{array}
\right):H_1\oplus H_2 \rightarrow  H_1\oplus H_2,B=\left(
\begin{array}{cc}B_1&B_{12}\\ 0&B_2\\ \end{array}
  \right):K_1\oplus K_2 \rightarrow  K_1\oplus K_2,$$

\noindent where $\dim (H_1)<\infty$, $\dim(K_2)<\infty,$ $A_2$ is
left invertible, $B_1$ is right invertible, and both $A_1$ and $B_2$
are nilpotent. Similar to the proof of Theorem 2.9, we have

\begin{equation} \alpha
(B_1)=\beta(A_2).\end{equation}
Moreover, considering the operator $C'\in B(K_1, H_2)$ with the
following expression

$$C' =\left(\begin{array}{cc} C'_1&C'_2\\
 C'_3&C'_4\\\end{array}\right): N(B_1)\oplus N(B_1)^\perp
\rightarrow R(A_2)^\perp \oplus R(A_2), $$ it follows from Corollary
2.11 (replacing operators $A, B, C$ by $A_1, B_1, C'$, respectively)
that
$$ \left(
       \begin{array}{cc}
        A_2&C'\\ 0&B_1\\
       \end{array}\right)
\makebox {is invertible if and only if so is}\,\, C'_1,$$ since
$A_2$ is left invertible and $B_1$ is right invertible.

\noindent Furthermore, since $\dim N(B_1)^\perp=\dim R(A_2)=\infty,$
we have
\begin{equation}\dim K_2+\dim N(B_1)^\perp=\dim R(A_2)+\dim
H_1=\infty.\end{equation} It follows from (5) and (6) that there
exist two invertible operators $Q_1\in B(N(B_1), R(A_2)^\perp)$ and
    $$Q_2=\left(
    \begin{array}{cc}
       T_1&T_2 \\
       T_3&T_4\\
    \end{array}\right):
    N(B_1)^\perp \oplus K_2
    \rightarrow R(A_2)\oplus H_1.$$
Define an operator

$$C : K\rightarrow H\, \makebox{ by}\, C = \left(
       \begin{array}{ccc}
       Q_1&0 &0\\
       0&T_1&T_2\\
       0& T_3&T_4\\
       \end{array}\right):
 N(B_1)\oplus N(B_1)^\perp \oplus K_2
 \rightarrow  R(A_2)^\perp\oplus R(A_2)\oplus H_1.$$
Obviously $C\in B(K,H)$ is invertible. Next we will prove that
$M_C$ is a Browder operator.
Observing that

$M_C=\left(\begin{array}{ccccc}
   A_1&A_{12}& 0&T_3&T_4\\
   0&0 & Q_1&0 &0\\
   0&A'_2 &0& T_1&T_2\\
   0&0 &0&B'_1&B_{12}\\
   0&0 &0& 0&B_2\\
  \end{array} \right)
:  H_1\oplus H_2 \oplus  N(B_1)\oplus N(B_1)^\perp \oplus K_2
  \rightarrow H_1\oplus R(A_2)^\perp\oplus R(A_2)\oplus  K_1\oplus K_2.$

  $=\left(\begin{array}{cccc}
 A_1& A_{12}&C_1&T_4\\
 0&A_2 & C_3&C_2\\
 0&0&B_1&B_{12}\\
 0&0& 0&B_2\
 \end{array}\right):
 H_1\oplus H_2 \oplus K_1\oplus K_2
 \rightarrow H_1\oplus H_2\oplus K_1\oplus K_2,$

\noindent where $A_2=\left(\begin{array}{c}
   0\\A'_2 \\ \end{array}\right),
   C_1=\left(\begin{array}{cc}
   0&T_3\\\end{array}\right),
   C_3=\left(\begin{array}{cc}
   Q_1&0 \\ 0&T_1\\
   \end{array}\right), C_4=\left(\begin{array}{c}
   0\\T_2 \\ \end{array}\right), B_1=\left(\begin{array}{cc}
   0&B'_1\\\end{array}\right).$
From Corollary 2.14 we have that $\left(\begin{array}{cc}
   A_2&C_3 \\ 0&B_1\\
   \end{array}\right)$ is invertibility, and then it follows from
   Lemma 2.4 that $M_C$ is a Browder operator. This ends the proof.

%
%
%

\end{document}